\documentclass[a4paper,twoside,12pt]{article}
\usepackage{amssymb,amsmath,amsthm,latexsym}
\usepackage{amsfonts}
\usepackage{graphicx,caption,subcaption}
\usepackage[pdftex,bookmarks,colorlinks=false]{hyperref}
\usepackage[hmargin=1.25in,vmargin=1.25in]{geometry}
\usepackage[mathscr]{euscript}
\usepackage[auth-sc-lg,affil-sl]{authblk}
\def\ni{\noindent}

\def\sP{\mathscr{P}}
\def\S{\mathscr{S}}
\def\X{\mathscr{X}}

\newtheorem{thm}{Theorem}[section]

\newtheorem{lem}[thm]{Lemma}

\title{\textbf{\sc A Note on Sparing Number Algorithm of Graphs}}

\author{N. K. Sudev}
\affil{\small Department of Mathematics\\ Vidya Academy of Science \& Technology \\  Thrissur - 680501, India.\\ E-mail: sudevnk@gmail.com}

\author{K. A. Germina}
\affil{\small Department of Mathematics\\ University of Botswana, Botswana.\\ E-mail: srgerminaka@gmail.com}

\date{}

\begin{document}
\maketitle

\begin{abstract}
Let $X$ denote a set of all non-negative integers and $\sP(X)$ be its power set. A weak integer additive set-labeling (WIASL) of a graph $G$ is an injective set-valued function $f:V(G)\to \sP(X)-\{\emptyset\}$ where induced function $f^+:E(G) \to \sP(X)-\{\emptyset\}$ is defined by $f^+ (uv) = f(u)+ f(v)$ such that either $|f^+ (uv)|=|f(u)|$ or $|f^+ (uv)|=|f(v)|$ , where $f(u)+f(v)$ is the sumset of $f(u)$ and $f(v)$. The sparing number of a WIASL-graph $G$ is the minimum required number of edges in $G$ having singleton set-labels.  In this paper, we discuss an algorithm for finding the sparing number of arbitrary graphs.
\end{abstract}

\vspace{0.2cm}

\ni \textbf{Key words}: Integer additive set-labeled graphs, weak integer additive set-labeled graphs, sparing number of a graph, spring number algorithm.

\vspace{0.04in}

\section{Introduction}

For all  terms and definitions, not defined specifically in this paper, we refer to \cite{BM1,FH, DBW}. Unless mentioned otherwise, all graphs considered here are simple, finite, non-trivial and connected.

The {\em sumset} of two sets $A$ and $B$ of integers, denoted by $A+B$, is defined as $A+B=\{a+b:a\in A, b\in B\}$.  If $A$ or $B$ is countably infinite, then their sumset $A+B$ will also be countably infinite. Hence, all sets we consider here are finite sets of non-negative integers.

Let $X$ be a non-empty finite set of non-negative integers and let $\sP(X)$ be its power set.  An {\em integer additive set-labeling} (IASL) of a graph $G$ (see \cite{GS1,GS0}) is an injective function $f:V(G)\to \sP(X)-\{\emptyset\}$ such that the induced function $f^+:E(G)\to \sP(X)-\{\emptyset\}$ is defined by $f^+{uv}=f(u)+f(v)~ \forall uv\in E(G)$.  A graph $G$ which admits an IASL is called an {\em integer additive set-labeled graph} (IASL-graph). 

The cardinality of the set-label of an element (vertex or edge) of a graph $G$ is called the {\em set-indexing number} of that element. An element of a given graph $G$ is said to be a {\em mono-indexed element} of $G$ if its set-indexing number is $1$.

A {\em weak integer additive set-labeling} of a graph $G$ is an IASL $f:V(G)\to \sP(X)-\{\emptyset\}$, where induced function $f^+:E(G) \to \sP(X)-\{\emptyset\}$ is defined by $f^+ (uv) = f(u)+ f(v)$ such that either $|f^+ (uv)|=|f(u)|$ or $|f^+ (uv)|=|f(v)|$ , where $f(u)+f(v)$ is the sumset of $f(u)$ and $f(v)$.

\begin{lem}\label{L-WIASLG1}
{\rm \cite{GS3}} An IASI $f:V(G)\to \sP(X)-\{\emptyset\}$ of a given graph $G$ is a weak IASI of $G$ if and only if at least one end vertex of every edge of $G$ mono-indexed.
\end{lem}

Hence, it can be seen that both end vertices of some edges of a given graph can be (must be) mono-indexed and hence those edges are also mono-indexed. The minimum number of mono-indexed edges required in a graph $G$ so that $G$ admits a WIASL is called the {\em sparing number} of $G$, denoted by $\varphi(G)$. 

Note that an independence set $I$ is said to have maximal incidence in $G$ if maximum number of edges in $G$ have their one end vertex in $I$. Then, the sparing number of any given graph can be determined using the following theorem.

\begin{thm}\label{T-WIASL-3}
{\rm \cite{GS2}} Let $G$ be a given WIASL-graph and $I$ be an independent set in $G$ which has the maximal incidence in $G$. Then, the sparing number of $G$ is the $|E(G-I)|$. 
\end{thm}

Certain studies on WIASL-graphs and their sparing numbers have been done in \cite{GS1,GS0,GS2,GS3,GS4}. In this paper, we discuss an algorithm to determine the sparing number of arbitrary finite connected graphs.

\section{Sparing Number Algorithm}

\ni In this section, we use the following notations.
%$X$ represents the ground set used for labeling the elements of the graph $G$. The neighbour set of a vertex $v$ in $G$ is denoted by $N(v)$ and $N[v]=\cup \{v\}$.

\begin{enumerate}\itemsep0mm
\item $X:=$ The ground set used for labeling the elements of the graph $G$.
\item $N(v):=$ The set of all vertices in $G$ adjacent to the vertex $v$.
\item $N[v]:= N(v)\cup \{v\}$.
\end{enumerate}

We now consider describe an algorithm to iteratively find out the sparing number of a given graph as explained below.

\subsection{The Sparing Number Algorithm}

\noindent\fbox{%
    \parbox{\textwidth}{
\begin{enumerate}\itemsep0.5mm

\item[(1)] Set $G_1=G, ~~ X_1=\emptyset,  ~~  Y_1=\emptyset,  ~~  E_1=\emptyset$.

\item[(2)] Choose $v_i$ such that $d(v_i)=\Delta(G_i)$. 

\item[(3)] Label the vertex $v_i$ by a non-empty, non-singleton subset of the ground set $X$. 

\item[(4)] Let $G_{i+1}=G_i-\{v_i\}$, $X_{i+1}=X_i\cup \{v_i\}$ and $Y_{i+1}=Y_i\cup N(v_i)$.

\item[(5)] For any two vertices $v_r, v_s \in N(v_i)$, if $v_rv_s\in E(G)$, then let $E_{i+1}=E_i\cup \{v_rv_s\}$. If no two vertices in $N(v_i)$ are mutually adjacent, then $E_{i+1}=E_i$. 

\item[(6)] Label the vertices in $Y_{i+1}$ by distinct singleton subsets of the ground set $X$.

\item[(7)] If all vertices of $G_i$ are labeled, then go to step (8). Otherwise, go to step (2).

\item[(8)] Here, the sparing number of $G$ is $\varphi(G)=|E_i|$. Stop.

\end{enumerate}
}}

\vspace{0.45cm}

Note that the set $I=\bigcup\limits_i X_i$ is the independent set of vertices in $G$ with maximal incidence in $G$ and hence by Theorem \ref{T-WIASL-3}, this algorithm provides the sparing number of any given graph.
Also, the vertices of $G$ in $I$ have non-singleton set-labels whereas the vertices in $V-I$ are mono-indexed. 

\subsection{An Illustrations to Sparing Number Algorithm} 

Consider the graph given below in Figure \ref{fig:Fig-1} for finding out the sparing number. Let $X$ be a set of non-negative integers, which is used as the ground set for labeling the vertices of $G$. Let $\S$ be the collection of all singleton subsets of $X$ and $\X$ be the collection of all non-empty, non-singleton subsets of the ground set $X$.  %We shall use the sparing number algorithm as explained below. 

\begin{figure}[h!]
\centering
\includegraphics[width=0.8\linewidth]{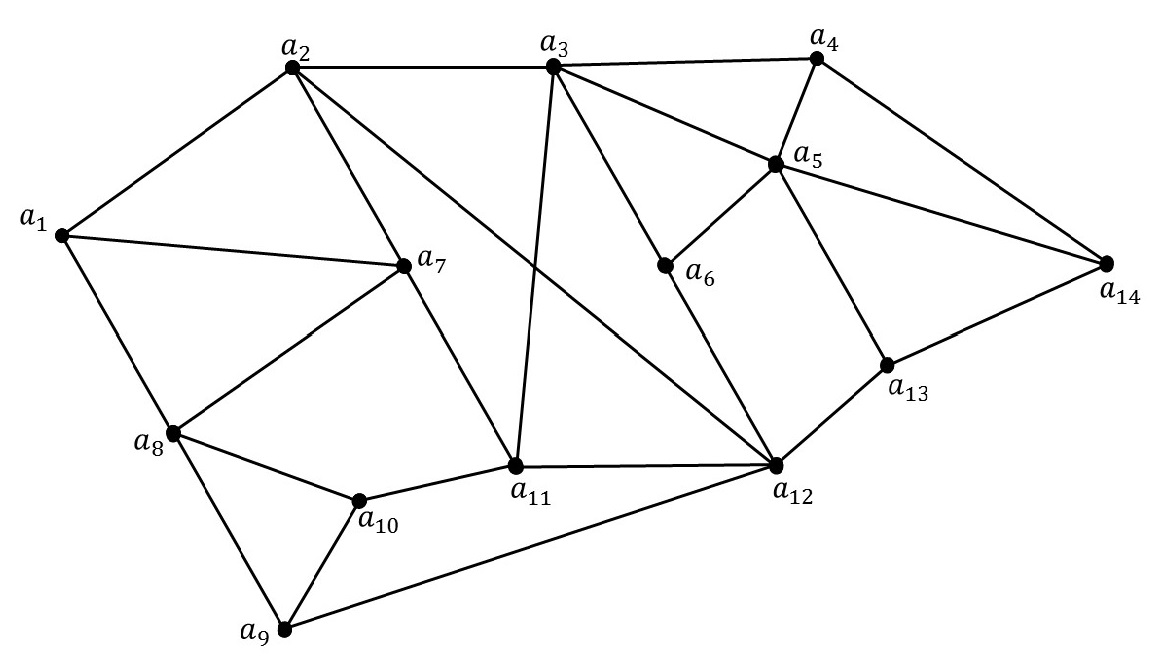}
\caption{}
\label{fig:Fig-1}
\end{figure}

%\begin{enumerate}
%\item[(i)] \textbf{First Iteration:}

First, let $G_1=G$, $X_1=\emptyset$, $Y_1=\emptyset$ and $E_1=\emptyset$. At this stage, we have $\Delta(G_1)=5$ and the vertices $a_3, a_5, a_{12}$ have degree $5$ in $G_1$.  Without loss of generality, let $v_1=a_3$ and $N(v_1)=\{a_2,a_4,a_5,a_6,a_{11}\}$.

Now, let $X_2=X_1\cup \{a_3\}= \{a_3\}$ and $Y_2=Y_1\cup N(v_1)=\{a_2,a_4,a_5,a_6,a_{11}\}$. Also, note that for the vertices, $a_4,a_5,a_6\in N(v_1)$, we have $a_4a_5, a_5a_6\in E(G)$ and hence $E_2=\{a_4a_5,a_5a_6\}$. 

Here, label the vertex $a_3\in X_2$ by a subset of $X$, that is in $\X$, and label the vertices of $Y_2$ by distinct singleton subsets of $X$ that are in $\S$.

Next, reduce the graph $G$ as $G_2$ such that $V(G_2)=V(G_1)-N[v_1]$. In this reduced graph, the vertex $a_{12}$ has the maximum degree, $d(a_{12})=5$. Hence, let $v_2=a_{12}$ and $N(v_2)=\{a_2,a_6,a_9,a_{11},a_{13}\}$.  Then, $X_3=X_2\cup \{a_{12}\}=\{a_3,a_{12}\}$ and $Y_3=Y_2\cup N(v_2)=\{a_2,a_4,a_5,a_6,a_9,a_{11},a_{13}\}$. Also, note that there is no two vertices in $N(v_2)$ are adjacent in $G_2$ and hence $E_3=E_2$. Now, label the vertex $a_{12}\in X_3$ by a subset of $X$ in $\X$ and label the unlabeled vertices in $Y_3$ (that is, the vertices in $Y_3-Y_2$), by distinct singleton subsets of $X$ in $\S$, which are not used for labeling in previous iterations.

Next, reduce the graph $G_2$ in to a new graph $G_3$ such that $V(G_3)=V(G_2)-N[v_2]$. In this reduced graph $G_3$, the vertex $a_7$ has the maximum degree, $d(a_7)=4$. Hence, let $v_3=a_7$ and then $N(v_3)= \{a_1,a_2,a_8,a_{11}\}$. Now, let $X_4=X_3\cup \{a_7\}=\{a_3,a_{12},a_7\}$ and $Y_4=Y_3\cup N(v_3)=\{a_1,a_2,a_4,a_5,a_6,a_9,a_{11},a_{13}\}$. Here, for the vertices, $a_1,a_2,a_8\in N(v_3)$, we have $a_1a_2, a_1a_8\in E(G)$ and hence $E_4=E_3\cup \{a_1a_2,a_1a_8\}=\{a_4a_5, a_5a_6,a_1a_2,a_1a_8\}$. Label the vertex $a_7\in X_4$ by a non-singleton set in $\X$, which is not used for labeling before and the vertices in $Y_4-Y_3$ by singleton sets in $\S$, which are not already used. 

Next, reduce the graph $G_3$ in to the graph $G_4$ such that $V(G_4)=V(G_3)-N[v_3]$.  In $G_4$, the vertex $a_{10}$ has the maximum degree $d(a_{10})=3$. Hence, let $v_4=a_{10}$ and then $N(v_{10})= \{a_8,a_9,a_{11}\}$. Now, let $X_5=X_4\cup \{a_{10}\}=\{a_3,a_{12},a_7,a_{10}\}$ and $Y_5=Y_4\cup N(v_3) =\{a_1,a_2,a_4,a_5,a_6,a_8,a_9,a_{11},a_{13}\}$. Here, for the vertices, $a_8,a_9\in N(v_4)$, we have $a_8a_9\in E(G)$ and hence $E_5=E_4\cup \{a_8a_9\}=\{a_4a_5, a_5a_6,a_1a_2,a_1a_8,a_8a_9\}$. 

Next, reduce the graph $G_4$ in to the graph $G_5$ such that $V(G_5)=V(G_4)-N[v_4]$. In $G_5$, the vertex $a_{14}$ has the maximum degree $d(a_{14})=3$. Hence, let $v_5=a_{14}$ and then $N(v_5)= \{a_4,a_5,a_{13}\}$. Now, let $X_5=X_4\cup \{a_{14}\}=\{a_3,a_{12},a_7,a_{10}\}$ and $Y_5=Y_4\cup N(v_4) =\{a_1,a_2,a_4,a_5,a_6,a_8,a_9,a_{11},a_{13}\}$. Here, for the vertices, $a_4,a_5,a_{13}\in N(v_5)$, we have $a_4a_5,a_5a_{13}\in E(G)$ and hence $E_5=E_4\cup \{a_4a_5,a_5a_{13}\}=\{a_4a_5,a_5a_6,a_1a_2,a_1a_8,a_8a_9,a_5a_{13}\}$. 

Next, reduce the graph $G_5$ in to the graph $G_6$ such that $V(G_6)=V(G_5)-N[v_5]$. This reduced graph is a trivial graph. Now, all the vertices of the given graph $G$ have been labeled by the subsets of $X$ in such a way that at least one end vertex of every edge in $G$ has singleton set-label. 

\begin{figure}[h!]
\centering
\includegraphics[width=0.85\linewidth]{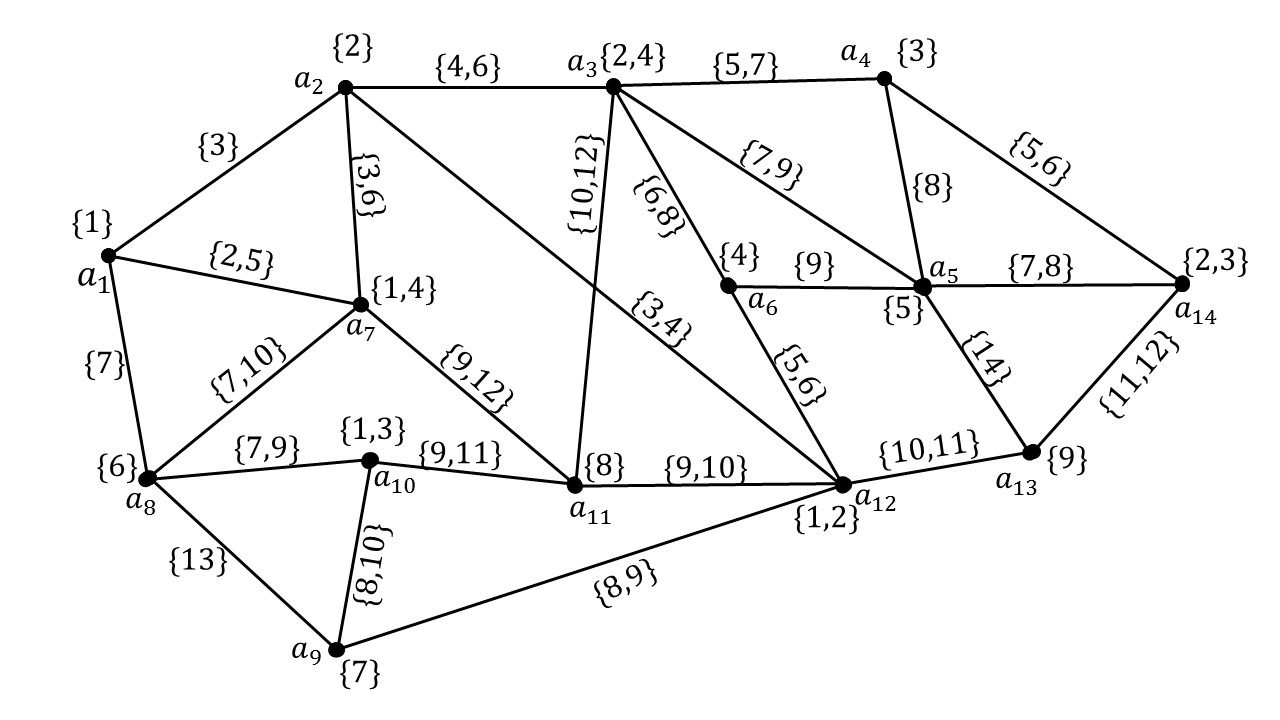}
\caption{}
\label{fig:Fig-2}
\end{figure}

Now, label the edges in $G$ in such a way that the set-label of every edge is the sumset of the set-labels of its end vertices. A labeling of the above graph $G$ as mentioned in the sparing number algorithm is illustrated in Figure \ref{fig:Fig-2}. 

Note that all edges listed in $E_5$ will be mono-indexed. Therefore, the sparing number of $G$ is given by $\varphi(G)=|E_5|=6$.

\section{Conclusion}

In this paper, we discussed an algorithm to determine the sparing number of arbitrary WIASL-graphs. A detailed study on the complexity of this algorithm seems to be possible and is interesting for further investigation. 

Further studies on many other characteristics of different IASL-graphs are also interesting and challenging. All these facts highlight the scope for further studies in this area.    

\section{Acknowledgements}

The authors would like to dedicate this paper to the bright memory of Prof. (Dr.) B. D. Acharya who introduced the concepts of set-valuations of graphs.


\begin{thebibliography}{20}

\bibitem{BDA1} B. D. Acharya, {\bf Set-valuations and their applications}, MRI Lecture Notes in Applied Mathematics, No.2, The Mehta Research Institute of Mathematics and Mathematical Physics, Allahabad, 1983.

\bibitem{ND} N. Deo, {\bf Graph theory with application to engineering and computer science}, Prentice Hall of India Pvt. Ltd., Delhi, 1974.

\bibitem{BM1} J. A. Bondy and U. S. R. Murty, {\bf Graph theory with applications}, North-Holland, New York, 1976.

\bibitem {GS1} K. A. Germina and N. K. Sudev,  {\em On weakly uniform integer additive set-indexers of graphs}, Int. Math. Forum, {\bf 8}(37)(2013), 1827-1834. DOI: 10.12988/imf.2013.310188.

\bibitem{FH}  F. Harary, {\bf Graph theory}, Addison-Wesley, 1969.

\bibitem {MBN} M. B. Nathanson,  {\bf Additive number theory, inverse problems and geometry of sumsets}, Springer, New York, 1996.

\bibitem{GS0} N. K. Sudev and K. A. Germina, {\em On integer additive set-indexers of graphs}, Int. J. Math. Sci. Engg. Appl., {\bf 8}(2)(2014), 11-22.

\bibitem{GS2} N. K. Sudev and K. A. Germina, {\em Some new results on weak integer additive set-labelings of graphs}, Int. J. Computer Appl., {\bf 128}(1)(2015),1-5., DOI: 10.5120/ijca2015906514.

\bibitem{GS3} N. K. Sudev and K. A. Germina, {\em A characterisation of weak integer additive set-indexers of graphs}, J. Fuzzy Set Valued Anal., {\bf 2014}(2014), 1-6., DOI:10.5899/2014/jfsva-00189

\bibitem{GS4} N. K. Sudev and K. A. Germina, {\em A note on the sparing number of graphs}, Adv. Appl. Discrete Math., {\bf 14}(1)(2014), 51-65.

\bibitem{GS5} N. K. Sudev and K. A. Germina, {\em Further studies on the sparing number of graphs}, TechS Vidya e-Journal of Research, {\bf 2}(1)(2014-15), 25-38.

\bibitem{DBW} D. B. West, {\bf Introduction to Graph Theory}, Pearson Education Inc., 2001.

\end{thebibliography}
\end{document}